\documentstyle{amsart}
\input{bull-art}
\makeatletter

\newif\if@defeqnsw \@defeqnswtrue

\def\eqnarray{\stepcounter{equation}\let\@currentlabel=%
\theequation
\if@defeqnsw\global\@eqnswtrue\else\global\@eqnswfalse\fi
\global\@eqnswtrue
\tabskip\@centering\let\\=\@eqncr
$$\halign to \displaywidth\bgroup\hfil\global\@eqcnt\z@
  $\displaystyle\tabskip\z@{##}$&\global\@eqcnt\@ne 
  \hfil$\displaystyle{{}##{}}$\hfil
  &\global\@eqcnt\tw@ $\displaystyle{##}$\hfil 
  \tabskip\@centering&\llap{##}\tabskip\z@\cr}

\def\yesnumber{\global\@eqnswtrue}

\def\@@eqncr{\let\@tempa\relax\global\advance\@eqcnt by \@ne
    \ifcase\@eqcnt \def\@tempa{& & & &}\or \def\@tempa{& & 
&}\or
     \def\@tempa{& &}\or \def\@tempa{&}\else\fi 
     \@tempa \if@eqnsw\@eqnnum\stepcounter{equation}\fi
     
\if@defeqnsw\global\@eqnswtrue\else\global\@eqnswfalse\fi
     \global\@eqcnt\z@\cr}

\def\@eqnacr{{\ifnum0=`}\fi\@ifstar{\@yeqnacr}{\@yeqnacr}}

\def\@yeqnacr{\@ifnextchar [{\@xeqnacr}{\@xeqnacr[\z@]}}

\def\@xeqnacr[#1]{\ifnum0=`{\fi}\cr 
\noalign{\vskip\jot\vskip #1\relax}}

\def\eqalign{\null\,\vcenter\bgroup\openup1\jot \m@th 
\let\\=\@eqnacr
\ialign\bgroup\strut
\hfil$\displaystyle{##}$&$\displaystyle{{}##}$\hfil\crcr}
\def\endeqalign{\crcr\egroup\egroup\,}

\def\cases{\left\{\,\vcenter\bgroup\normalbaselines\m@th 
\let\\=\@eqnacr
    \ialign\bgroup$##\hfil$&\quad##\hfil\crcr}
\def\endcases{\crcr\egroup\egroup\right.}

\def\eqalignno{\stepcounter{equation}\let\@currentlabel=%
\theequation
\if@defeqnsw\global\@eqnswtrue\else\global\@eqnswfalse\fi
\let\\=\@eqncr
$$\displ@y \tabskip\@centering \halign to 
\displaywidth\bgroup
  \global\@eqcnt\@ne\hfil
  $\@lign\displaystyle{##}$\tabskip\z@skip&\global\@eqcnt%
\tw@
  $\@lign\displaystyle{{}##}$\hfil\tabskip\@centering&
  \llap{\@lign##}\tabskip\z@skip\crcr}

\def\endeqalignno{\@@eqncr\egroup
      \global\advance\c@equation\m@ne$$\global\@ignoretrue}

\@namedef{eqalignno*}{\@defeqnswfalse\eqalignno}
\@namedef{endeqalignno*}{\endeqalignno}

\def\eqaligntwo{\stepcounter{equation}\let\@currentlabel=%
\theequation
\if@defeqnsw\global\@eqnswtrue\else\global\@eqnswfalse\fi
\let\\=\@eqncr
$$\displ@y \tabskip\@centering \halign to 
\displaywidth\bgroup
  \global\@eqcnt\m@ne\hfil
  $\@lign\displaystyle{##}$\tabskip\z@skip&\global\@eqcnt\z@
  $\@lign\displaystyle{{}##}$\hfil\qquad&\global\@eqcnt\@ne
  \hfil$\@lign\displaystyle{##}$&\global\@eqcnt\tw@
  $\@lign\displaystyle{{}##}$\hfil\tabskip\@centering&
  \llap{\@lign##}\tabskip\z@skip\crcr}

\def\endeqaligntwo{\@@eqncr\egroup
      \global\advance\c@equation\m@ne$$\global\@ignoretrue}

\@namedef{eqaligntwo*}{\@defeqnswfalse\eqaligntwo}
\@namedef{endeqaligntwo*}{\endeqaligntwo}

\newtoks\@stequation

\def\subequations{\refstepcounter{equation}%
  \edef\@savedequation{\the\c@equation}%
  \@stequation=\expandafter{\theequation}
  \edef\@savedtheequation{\the\@stequation}
  \edef\oldtheequation{\theequation}%
  \setcounter{equation}{0}%
  \def\theequation{\oldtheequation\alph{equation}}}

\def\endsubequations{%
  \setcounter{equation}{\@savedequation}%
  \@stequation=\expandafter{\@savedtheequation}%
  \edef\theequation{\the\@stequation}%
  \global\@ignoretrue}

\def\big#1{{\hbox{$\left#1\vcenter 
to1.428\ht\strutbox{}\right.\n@space$}}}
\def\Big#1{{\hbox{$\left#1\vcenter 
to2.142\ht\strutbox{}\right.\n@space$}}}
\def\bigg#1{{\hbox{$\left#1\vcenter 
to2.857\ht\strutbox{}\right.\n@space$}}}
\def\Bigg#1{{\hbox{$\left#1\vcenter 
to3.571\ht\strutbox{}\right.\n@space$}}}

\makeatother

\def\theequation{\thesection.\arabic{equation}}

\def\bbbld{\Bbb}
\newcommand{\bfR}{\bbbld R}
\newcommand{\bfC}{\bbbld C}

\def\Re{\operatorname{Re}}  

\newcommand{\n}{\noindent}
\newcommand{\p}{\hspace*{\parindent}}

\newtheorem{oldtheorem}{Theorem}
\newtheorem{oldassertion}[oldtheorem]{Assertion}
\newtheorem{oldproposition}[oldtheorem]{Proposition}
\newtheorem{oldremark}[oldtheorem]{Remark}
\newtheorem{oldlemma}[oldtheorem]{Lemma}
\newtheorem{olddefinition}[oldtheorem]{Definition}
\newtheorem{oldclaim}[oldtheorem]{Claim}
\newtheorem{oldcorollary}[oldtheorem]{Corollary}

\newbox\qedbox
\setbox\qedbox=\hbox{$\Box$}
			{\hfill\penalty10000\copy\qedbox\par\medskip}
		      {\smallskip\noindent{\bf Proof of #1.}\hskip 
\labelsep}%
			{\nobreak\hfill\hfill\nobreak\copy\qedbox\par\medskip}

\begin{document}
\def\currentvolume{29}
\def\currentissue{1}
\def\currentyear{1993}
\def\currentmonth{July}
\def\copyrightyear{1993}
\def\currentpages{77-84}

\title{Adding Handles to the Helicoid}



\author{David Hoffman}
\address{Department of Mathematics \\ University of 
Massachusetts \\
Amherst, Massachusetts 01003}
\email{david@@gang.umass.edu}
\thanks{The first author was supported by the National 
Science
Foundation under 
grants, DMS-9101903 and DMS-9011083 and by the U.S. 
Department of
Energy under grant DE-FG02-86ER25015 
of the Applied Mathematical Science subprogram of the 
Office of Energy
Research. The second author was supported by the U.S. 
Department of
Energy under grant DE-FG02-86ER25015 
of the Applied Mathematical Science subprogram of the 
Office of Energy
Research. The third author was partially supported by 
Sonderforschungsbereich
SFB256 at Bonn} 

\author{Fusheng Wei}
\address{Department of Mathematics \\ Texas A\&M 
University\\
College Station, Texas 77843}
\curraddr{Department of Mathematics, University of 
Massachusetts, Amherst, 
Massachusetts
01003}
\email{wei@@gang.umass.edu}
\author{Hermann Karcher}
\address{Mathematisches Institut \\ Universitat Bonn \\
D-5300 Bonn, West Germany} 
\email{UNM416@@IBM.rhrz.uni-bonn.de}
\date{July 1, 1992 and, in revised form, December 6, 1992}

\maketitle

\begin{abstract}
	There exist two new embedded minimal surfaces, asymptotic 
to the
helicoid. One is periodic, with quotient (by
orientation-preserving translations) of genus one.  The 
other is nonperiodic 
of genus one.
\end{abstract}

\p	We have constructed two minimal surfaces of theoretical
interest.  The first is a complete, embedded, singly 
periodic minimal 
surface (SPEMS) that is asymptotic to the helicoid, has 
infinite genus,  
and whose quotient by translations has genus one. The 
quotient of 
the helicoid by translations has genus zero and the 
helicoid itself is 
simply connected. 
\vspace{.13in}

\n {\bf Theorem~1.}\, 
{\em There exists an embedded singly periodic minimal 
surface  
$ {\cal W}_{1} $, asymptotic to the helicoid and invariant 
under a 
translation $T$. The quotient surface $ {\cal W}_{1}/T $ 
has genus
equal to one 
and two ends.}
\vspace{.13in}

$ {\cal W}_{1} $ contains a vertical axis, as does the 
helicoid, and 
$ {\cal W}_{1}/T $ contains two horizontal lines.
The second surface is a complete, properly embedded minimal
surface of finite topology with infinite total curvature.  
It is the
first such surface to be found since the helicoid, which 
was discovered
in the eighteenth century.   (See Figure 1(a).)
\vspace{.13in}

\n {\bf Theorem~2.}\,
{\em There exists a complete, properly embedded minimal 
surface, 
$ {\cal H}\!e_1 $, of genus-one, whose one end is of 
helicoidal type.}
\vspace{.13in}

$ {\cal H}\!e_1 $ contains a vertical line, like the 
helicoid, and one
horizontal line that crosses it.  Schwarz reflection in 
these two lines 
generates the symmetry group of the surface. (See  Figure 
1(b).)

\begin{figure}[t]
\vspace{23.5pc}
\caption{The surfaces: (a) $ {\cal W}_{1} $ and (b) $ 
{\cal H}\!e_1 $.}
\end{figure}

\section{History and context}
\setcounter{equation}{0}
\p	Except for the plane, the helicoid is the only ruled 
minimal surface. 
Its discovery is attributed to Meusnier in 1776; together 
with the 
catenoid (Euler, c.\ 1744) these were the only minimal 
surfaces 
explicitly known to eighteenth century mathematics.  
The next major discovery came from Scherk in  
the 1830s: multiple families of 
periodic minimal surfaces, including the famous families of 
singly and doubly periodic examples that bear his name
\cite{ni2,sche1}.  That the surfaces  
in these two families share the same normal mapping is 
implicit in 
the work of Scherk. This fundamental relationship was made 
explicit 
by Enneper, Weierstrass, and Riemann. They developed an 
integral 
representation formula for minimal surfaces via contour 
integration 
of meromorphic data derived from the normal mapping, which 
they 
knew to be conformal \cite{os1}.  (See (3.1) below.)  
Minimal 
surfaces were seen to be, from this  
point of view, the real part of  complex curves in $ 
\bfC^3 $. The helicoid and 
catenoid were recognized as, locally, the real and 
imaginary parts of 
the same curve. The same is true of the two families of 
Scherk. 
Minimal surfaces related in this way are said to be {\em 
conjugate}.

	For complete minimal surfaces, whose quotient by 
orientation-preserving symmetries has finite total 
curvature, the
quotient is naturally a compact Riemann surface, possibly 
punctured
in a finite number of points.  Moreover, the meromorphic 
data is
well defined on the compact surface (Osserman 
\cite{os3,os1}).
Translations are produced when the integral representation 
has
periods on the Riemann surface.
The classical examples mentioned above can be represented 
on a 
sphere punctured two (resp. four) times for the 
helicoid/catenoid 
pair (resp. Scherk's singly/doubly periodic surfaces).  
Moreover, the 
Gauss map can be taken to be the identity on $ S^2 $. All 
these examples 
are {\em embedded}. 

The existence of higher-genus embedded examples has been an 
open question until recently. For complete embedded 
examples of finite
total curvature (FEMS) in $ {\bfR}^{3} $, Lopez-Ros 
\cite{lor1} showed
that the plane and the catenoid are the only FEMS of genus 
zero.
Schoen \cite{sc1} proved  that the catenoid was the only 
FEMS with two 
ends. The 
existence of examples with genus greater than zero and 
more than 
two ends is well documented [2--5, 18].

	The helicoid was the only known example of a complete 
embedded minimal surface with finite topology and infinite 
total 
curvature.  It has been a longstanding open question as to 
whether there 
are others.  Theorem~2 answers this question 
affirmatively.  Also, 
all known properly embedded examples of infinite
total curvature had infinite symmetry groups, and the 
quotients of
these surfaces by these groups were compact (possibly 
punctured)
Riemann surfaces whose inherited metric had finite total 
curvature.
The surface $ {\cal H}\!e_{1} $ of Theorem~2 is 
conformally a 
once-punctured rhombic torus with symmetry group $ 
Z_{2}\oplus
Z_2 $.
 
All classical, complete, embedded, doubly periodic minimal 
surfaces  
(DPEMS) can be defined by meromorphic data with periods on 
punctured spheres.
Karcher \cite{ka7} and Meeks-Rosenberg \cite{mr3} 
constructed new
families of DPEMS that had genus one in the quotient.  No 
higher-genus examples were known that were not coverings 
of these
examples. Moreover, there were no known  
genus-one examples with the same end behavior as the Scherk 
doubly periodic examples. In \cite{wei2} Wei constructed 
the first DPEMS of genus equal to two in the quotient.  
Based on 
the construction 
strategy used in that paper, Karcher was able to modify 
Scherk's 
doubly periodic example to produce a genus-one DPEMS that 
had 
the same end behavior as the Scherk example.  We refer to 
this 
surface as $ {\cal K}_{{\pi}/{2}} $ for reasons that will 
be made 
clear.   (See Figure 2.)

\begin{figure}[t]

\vspace*{30.5pc}
\caption{Scherk's doubly periodic surface: (a)
$\theta = {{\pi}/{2}}$ and  (b)
${\cal K}_{\frac{\pi}{2}}$.}
\end{figure}

\section{SPEMS as limits of DPEMS}
\setcounter{equation}{0}
\p	The Scherk family can be considered to be the 
desingularization 
of two 
families of equally spaced, parallel, vertical halfplanes 
meeting at an 
angle $\theta, 0<\theta \leq \pi/2$. In the slab 
$|x_{3}|<\epsilon$, the 
surfaces look like saddles over alternating regions in a 
tiling of $x_{3}
= 0$ by rhombi. (See Figures 2(a) and 3(a).)
With appropriate scaling as $\theta$ goes to zero, the 
rhombi diagonals grow 
in one direction only, approaching a strip in the plane.  
There is a basic 
relationship between the Scherk family and the helicoid. 
Namely, if 
one keeps the symmetric point of a fixed saddle at the 
origin, the 
limit surface, with appropriate scaling as $\theta$ goes 
to zero, exists 
and is the helicoid (Hoffman and Wohlgemuth \cite{hw1}).   

The generalization, $ {\cal K}_{{\pi}/{2}} $, of Scherk's 
surface  can be
understood as  Scherk's surface with a tunnel replacing 
every other saddle.
(See Figure 2(b).) The  underlying Riemann surface is the 
square torus punctured
in four  points. We proved that this surface can be 
deformed  in exactly the 
same manner as the Scherk family.  (See Figure 3 on page 
81.)
\vspace{.13in}

\begin{figure}[t] 

\vspace*{36pc}
\caption{A surface in the Scherk family: (a)  
${\cal K}_{\theta}$  and (b)
$\theta = {\pi}/{4}$.}
\end{figure}

\n {\bf Proposition~1.}\, 
{\em There exists a one-parameter family $ {\cal 
K}_{\theta}$ 
of embedded doubly periodic minimal surfaces, whose 
quotient has 
genus equal to one and four Scherk ends, two up and two 
down.  Each  
genus-one surface is a rhombic torus. The angle $\theta$ 
between the up 
and down ends, $ 0<\theta\leq \pi/2 $, parametrizes the 
family.}
\vspace{.13in}

Each member of the family may be considered to be a 
desingularization of two families of parallel halfplanes. 
Unlike the 
Scherk family,  these planes are not equally spaced. The 
interplanar 
distance alternates. The smaller distance between planes 
is spanned 
by tubes, while the larger one is bridged by saddles.

The singly periodic surface $ {\cal W}_1 $ of Theorem~1 
has the same 
relationship to $ {\cal K}_\theta $ as the helicoid has to 
the 
Scherk family.  Namely, 
choose a distinguished point in a fundamental domain that 
is 
identifiable on each surface (e.g., the point on the 
saddle  
where the normal is vertical) and keep this point at the 
origin.  Then 
\vspace{.14in}

\n {\bf Theorem~3.}
{\em The limit surface  as $\theta \rightarrow 0$ of the 
surfaces 
${\cal K}_\theta$ exists and is equal to $ {\cal W}_1 $.}

\section{Construction of $ {\cal H}\!e_{1} $}
\setcounter{equation}{0}
\p	The surface $ {\cal W}_{1} $ can be described as a 
helicoid,
into which has been sewn a handle at every other 
half-turn.  Thus a 
handle has been added to the surface modulo translation.  
One 
could imagine adding a handle to every other fundamental 
domain,
producing three half-twists between handles, and (why stop 
at 
three?) in general $ 2k+1 $ half-twists between handles, $ 
k\geq 0 $.
The quotient by orientation-preserving translations of 
such a 
surface will have genus-one.  Now imagine fixing one 
horizontal
line in a fundamental domain to be the $ x_{2}$-axis and 
letting
$ k\rightarrow \infty $.  The resulting surface will have 
genus
one, will contain the $ x_{2}$-axis and $ x_{3}$-axis but 
no other
lines, will {\em not} be periodic, and should be 
asymptotic in some
sense to the helicoid.  In fact, such a surface exists and 
is the
surface $ {\cal H}\!e_{1} $ of Theorem~2.  Because of the 
presence of two
lines on $ {\cal H}\!e_{1} $ that cross at one point on 
the surface (and
at the end), $ {\cal H}\!e_{1} $ is invariant under two 
anticonformal reflections
with precisely two common fixed points.  This implies that 
$ {\cal H}\!{e_{1}} $
is conformally a rhombic torus and can be written as
\[w^{2}=(z-\lambda) (z^{2}+1)\]
for some $ \lambda \in {\bfR} $, with the end at $ 
(w,z)=(\infty,
\infty) $.  Note that $ {dz}/{w} $ has no zeros or poles.

	The key to constructing Weierstrass data on $ {\cal 
H}\!e_{1} $ is to 
realize 
that, while its Gauss map $ g $ has an essential 
singularity at the
end, its logarithmic differential, $ {dg}/{g} $, is 
meromorphic.
Note that, in general, at an end where $ g $ has a pole or 
a zero of finite
order, $ {dg}/{g} $ has a simple pole.  The helicoid may be
described on $ {\bfC} $ by the Weierstrass data $ g=e^{z}, 
dh=i dz $.  The single end occurs at infinity, where $ 
{dg}/{g} $
has a {\em double} pole, as does $ dh $.  Thus, in order 
to expect a
helicoidal end on a torus, we must look for a meromorphic 
differential 
$ dh $ with a double pole at the end and, by the Riemann 
relation, two 
zeros.  In order to have a helicoidal end at $ (\infty, 
\infty) $, we
specify that $ {dg}/{g} $ and $ dh $ should have double 
poles there,
as this is what happens on the helicoid.  The Riemann 
relation, the 
expected symmetries of the surface, the necessity for $ 
{dg}/{g} $ to
have integer residues and periods that are integer 
multiples of $ 2\pi $
(because $ g=e^{\int {dg}/{g}} $ must be well defined on $ 
{\cal H}\!
e_{1} $), together with the aforementioned 
Gedanken-experiment, allow us to
place the zeros and poles of $ {dg}/{g} $ (which are 
branch points and
vertical points of $ g )$.  We get
\[ 
\frac{dg}{g}=\rho \frac{(z-\alpha) 
(z-\beta)}{(z-a)}\,\frac{dz}{w} , \qquad 
dh=i(z-a)\,\frac{dz}{w}\,.\]
Here $ a,\alpha $, and $ \beta $ are real, and $ 
\rho=\left({(a-\alpha)
(a-\beta)}/{w(a)}\right)^{-1} $.  The two points on $ 
{\cal H}\!e_{1} $
where $ g $ is vertical are the zeros of $ dh $ and occur 
when $ z=a $,
where $ {dg}/{g} $ has simple poles.  The points on $ 
{\cal H}\!e_{1} $ 
where $ z=\alpha,\beta $ are branch points of $ g $ and 
zeros of $ {dg}/
{g} $.  Period conditions on $ {dg}/{g} $ imply that $ 
\alpha $ and
$ \beta $ are solutions to a quadratic expression with 
coefficients that
are functions of $ \lambda $ and $ a $.  Period conditions 
for the
Weierstrass integral
\begin{equation}
X(p)=\Re\int^{p}_{p_{0}}(\tfrac{1}{2}(g^{-1}-g) \,dh, 
\tfrac{i}{2} (g^{-1}
+g) dh, dh)
\end{equation}
to be well defined on $ {\cal H}\!e_{1} $ determine $ a $ 
as a function of
$ \lambda $ and give one remaining real-valued function of 
$ \lambda $ that
must vanish.  For $ \lambda \sim 0.32 $ this condition is 
satisfied.

\section{Computation}
\setcounter{equation}{0}
\p	To produce the pictures and find strong experimental 
evidence
that  these surfaces exist, MESH \cite{chh1,h2} was used.
Computational  programs were used to solve the 
period problems inherent in these representations.  
A proof for the existence of $ {\cal W}_1 $ can be 
given using a degree theory argument for the period 
mapping. 
\vspace{.2in}

\section{Acknowledgments} 
\p The authors thank Harold Rosenberg for encouraging 
us to look for surfaces like $ {\cal W}_{1} $ and $ {\cal 
H}\!e_{1} $.  
The program MESH developed 
at GANG by Jim Hoffman et al.\ was essential to our work; 
without it 
we would not have been able to discover this surface $ 
{\cal W}_{1} $
by our methods.  Hoffman prepared the illustrations for 
this paper.
We were not the first to think about the existence of a 
minimal surface like $ {\cal W}_1 $. In 1985,  after 
viewing  images 
of the helicoid and 
Costa's surface in a news article \cite{pet}, Sarah Mellon 
constructed 
models of more complicated embedded minimal surfaces of 
finite 
total curvature out of wire and strapping tape, and a 
higher genus 
helicoid out of clay.  She got the finite total curvature  
surfaces more 
or less exactly right. The generalized helicoid she 
constructed was not  
${\cal W}_{1}$, but the morphology was correct.  

%

\end{document}